\documentclass{ifacconf}

\usepackage{graphicx}      
\usepackage{natbib}        
\usepackage{amsmath, amssymb}
\usepackage{setspace}
\usepackage{tikz}
\usetikzlibrary{matrix,arrows}
\usetikzlibrary{shapes}
\usetikzlibrary{calc}
\usetikzlibrary{trees,snakes}
\usetikzlibrary{positioning} 
\usetikzlibrary{decorations.markings}
\usetikzlibrary{patterns}
\definecolor{CornflowerBlue}{rgb}{0.258824,0.258824,0.435294}
\definecolor{SkyBlue}{rgb}{0.196078,0.6,0.8}
\definecolor{dblue}{rgb}{.098,.243,.424}
\definecolor{lblue}{rgb}{.33,.57,.835}
\definecolor{llblue}{rgb}{.447,.643,.831}
\definecolor{mblue}{rgb}{0.176, 0.380, 0.659}
\definecolor{lcomp}{rgb}{.969,.765,.416}
\definecolor{ddorange}{rgb}{0.624, 0.365, 0}
\definecolor{dorange}{rgb}{0.72, 0.506, 0.125}
\definecolor{lorange}{rgb}{0.961, 0.678, 0.165}
\definecolor{lgreen}{rgb}{.812,.969,.435}
\definecolor{lyellow}{rgb}{1,.859,.451}
\definecolor{dyellow}{rgb}{.651,.482,0}
\definecolor{lrred}{rgb}{1,.6,.451}
\definecolor{dred}{rgb}{.85,.176,0}
\colorlet{mmblue}{blue!50}
\newcommand{\lightercolor}[3]{
    \colorlet{#3}{#1!#2!white}
}
\newcommand{\darkercolor}[3]{
    \colorlet{#3}{#1!#2!black}
}
\lightercolor{mmblue}{50}{llblue}
\darkercolor{mmblue}{50}{ddblue}
  \colorlet{dgreen}{green!40!black}
\tikzset{
    >=stealth',
    punkt/.style={
           rectangle,
           rounded corners,
           draw=black, very thick,
           text width=6.5em,
           minimum height=2em,
           text centered},
    pil/.style={
           ->,
           thick,
           shorten <=2pt,
           shorten >=2pt,}
}

\tikzset{
  big arrow/.style={
    decoration={markings,mark=at position 1 with {\arrow[scale=0.5, #1]{>}}},
    postaction={decorate},
    shorten >=0.4pt, line width=0.5mm},
  }

\newcommand{\vep}{\varepsilon}
\newcommand{\rar}{\rightarrow}
\newcommand{\mb}{\mathbb}
\newcommand{\mc}{\mathcal}
\newtheorem{problem}{Problem}
\newcommand{\bmat}[1]{\begin{bmatrix}#1\end{bmatrix}}

\newcommand{\yref}{y^{\mathrm{ref}}}
\newcommand{\y}[1]{y^{(#1)}}
\newtheorem{proposition}{Proposition}

\newtheorem{remark}{Remark}
\DeclareMathOperator*{\argmax}{arg\,max}
\newtheorem{example}{Example}
\newcommand{\vmax}{\bar{v}}
\newcommand{\ymax}{\bar{y}}
\begin{document}
\begin{frontmatter}

  \title{Incentive Design and Utility Learning via Energy Disaggregation}

\thanks[footnoteinfo]{The work presented is supported by the NSF
Graduate Research Fellowship under grant DGE 1106400, NSF
CPS:Large:ActionWebs award number 0931843, TRUST (Team for Research in
Ubiquitous Secure Technology) which receives support from NSF (award
number CCF-0424422), and FORCES (Foundations Of Resilient
CybEr-physical Systems), the European Research Council
   under the advanced grant LEARN, contract 267381, a postdoctoral grant from the Sweden-America
   Foundation, donated by ASEA's Fellowship Fund, and  by a postdoctoral
   grant from the Swedish Research Council.
}

\author[All]{Lillian J. Ratliff} 
\author[All]{Roy Dong} 
\author[All]{Henrik Ohlsson}
\author[All]{S. Shankar Sastry}
\address[All]{University of California, Berkeley, 
Berkeley, CA 94720 USA {\tt\{ratliffl, roydong, ohlsson, sastry\}@eecs.berkeley.edu}}

\begin{abstract}                
The utility company has many motivations for modifying energy consumption patterns of consumers such as revenue decoupling and demand response programs. We model the utility company--consumer interaction as a reverse Stackelberg game and present an iterative algorithm to design incentives for consumers while estimating their utility functions. Incentives are designed using the aggregated as well as the disaggregated (device level) consumption data. We simulate the iterative control (incentive design) and estimation (utility learning and disaggregation) process for examples.
\end{abstract}

\begin{keyword}
Game Theory, Economic Design, Energy Management Systems
\end{keyword}

\end{frontmatter}

\section{Introduction}
\label{sec:intro}

Currently, most electricity distribution systems only provide aggregate power consumption feedback to consumers, in the form of a energy bill. Studies have shown that providing device-level feedback on power consumption patterns to energy users can modify behavior and improve energy efficiency~\citep{gardner:2008aa,laitner:2009aa}. 

The current infrastructure only has sensors to measure the aggregated power consumption signal for a household. Even advanced metering infrastructures currently being deployed have the same restriction, albeit at high resolution and frequency~\citep{armel:2013aa}. Additionally, deploying plug-level sensors would require entering households to install these devices. Methods requiring plug-level sensors are often referred to as \emph{intrusive load monitoring}, and the network infrastructure required to transmit high resolution, high frequency data for several devices per household would be very costly.

 A low cost alternative to the deployment of a large number of sensors is \emph{non--intrusive load monitoring}. We consider the problem of nonintrusive load monitoring, which, in the scope of this paper, refers to recovering the power consumption signals of individual devices from the aggregate power consumption signal available to our sensors. This is also sometimes referred to as \emph{energy disaggregation}, and we will use the two terms interchangeably. This problem has been an active topic of research lately. Some works include \cite{dong:2013ac,froehlich:2011aa,johnson:2012aa}.

We propose that the utility company should use incentives to motivate a change in the energy consumption of consumers. 
We assume the utility company cares about the satisfaction of its consumers as well as altering consumption patterns, but it may not be able to directly observe the consumption patterns of individual devices or a consumer's satisfaction function.

In brief, the problem of behavior modification in energy consumption can be understood as follows. The utility company provides incentives to myopic energy consumers, who seek to maximize their own utility by selecting energy consumption patterns. This can be thought of as a control problem for the utility company. Additionally, the utility company does not directly observe the energy consumption patterns of individual devices, and seeks to recover it from an aggregate signal using energy disaggregation. This can be thought of as an estimation problem. Further, the consumer does not report any measure of its satisfaction directly to the utility. Thus, it must be estimated as well.

There are many motivations for changing energy consumption patterns of users. Many regions are beginning to implement revenue decoupling policies, whereby utility companies are economically motivated to decrease energy consumption~\citep{eom:2008aa}.  Additionally, the cost of producing energy depends on many variables, and being able to control demand can help alleviate the costs of inaccurate load forecasting. Demand response programs achieve this by controlling a portion of the demand at both peak and off-peak hours~\citep{mathieu:2012aa}. We propose a model for how utility companies would design incentives to induce the desired consumer behavior.


%
In this paper, we consider three cases of utility learning and incentive design. In the first, the utility company designs an incentive based entirely on the aggregate power consumption signal. 
We propose an algorithm to estimate the satisfaction function of the consumer based on the consumer's aggregated power consumption signals in Section~\ref{sec:aggregate}.
Then, in Section~\ref{sec:disaggregation_a}, we consider the case where the utility company knows the power consumption signal of individual devices and an unknown satisfaction function. 
Finally, in Section~\ref{sec:disaggregation_b}, we consider the case when the utility company only has access to the aggregated power consumption signal, and uses an energy disaggregation algorithm to recover the power consumption of individual devices. This disaggregated signal is used to allocate incentives, but the results will depend on the accuracy of our estimator, the energy disaggregation algorithm. 
We conclude the paper by showing the results from simulations of two examples of designing incentives while estimating the consumer's satisfaction function in Section~\ref{sec:numerical}. 
   Finally, in Section~\ref{sec:discussion} we make concluding remarks and discuss future research directions.


\begin{figure}[ht]
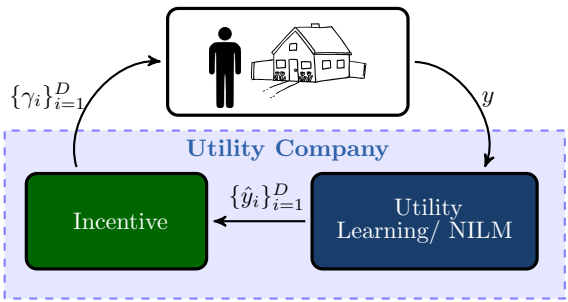

  \begin{center}
    \include{figs/loop}
  \end{center}
  \caption{Closing the Loop: Behavior modification via incentives $\gamma_i$ is a \emph{control} problem. The consumer decides when to use devices resulting in device level consumption $y_i$. Non--intrusive load monitoring (NILM) is used to \emph{estimate} problem device level usage $\hat{y}_i$. Similarly, utility learning is an estimation problem}
  \label{fig:loop}
\end{figure}

\section{Incentive Design Preliminaries}
A reverse Stackelberg game is a hierarchical control problem in which sequential decision making occurs; in particular, there is a leader that announces a mapping of the follower's decision space into the leader's decision space, after which the follower determines his optimal decision variables~\citep{groot:2012aa}.

Both the leader and the follower wish to maximize their pay--off determined by the functions $J_L(v, y)$ and $J_F(v,y)$ respectively. The leader's decision is denoted $v$; the follower's decision, $y$; and the incentive, $\gamma:y\mapsto v$. The basic approach to solving the reversed Stackelberg game is as follows. Let $v$ and $y$ take values in $V\subset \mb{R}$ and $Y\subset \mb{R}$, respectively; $J_L: \mb{R}\times \mb{R}\rar \mb{R}$; $J_F:\mb{R}\times \mb{R}\rar \mb{R}$. We define the desired choice for the leader as
\begin{equation}
  (v^d,y^d)=\arg\max_{v,y}\big\{J_L(v,y)| \ v\in V, y\in Y\}.
  \label{eq:desired}
\end{equation}
The incentive problem can be stated as follows:
\begin{problem}
  \label{prob:incentive}
  Find $\gamma:Y\rar V$, $\gamma\in \Gamma$ such that
  \begin{align}
    \arg\max_{y\in Y} J_F(\gamma(y), y)&= y^d\\
    \gamma(y^d)&=v^d
    \label{eq:prob1}
  \end{align}
  where $\Gamma$ is the set of admissible incentive mechanisms.
  \end{problem}

\section{Incentive Design Using Aggregate Power Signal}
\label{sec:aggregate}
We cast the utility--consumer interaction in a reversed Stackelberg game framework in which the utility company is the leader and the consumer is the follower (see Figure~\ref{fig:loop}). 
The leader's true utility is assumed to be given by
\begin{align}
  J_L(v,y)=g(y)-v+\beta f(y)
  \label{eq:trueprincip}
\end{align}
where $g(\cdot)$ is a concave function of the consumer's energy usage $y$ over a billing period, $v$ is the value of the incentive paid to the consumer, $f:Y\rar\mb{R}$ is the consumer's satisfaction function for energy consumption which we assume is concave and $\beta$ is a multiplying factor capturing the degree of benevolence of the utility company. 

We assume that $v\in V=[0,\vmax]$ since the utility company should not take additional money away from the consumer on top of the cost of their usage and $v$ should be less than some maximal amount the leader is willing to pay to the consumer $\vmax$. 
Similarly, let $y\in Y=[0, \ymax]$ where $\ymax$ is the upper bound on the allowed energy usage and let $\mathring{Y}=(0, \ymax)$. 

In a regulated market with revenue decoupling in place, a simplified model may consider 
\begin{equation}
  g(y)=-y
  \label{eq:revdecoup}
\end{equation}
 representing the fact that the utility wants the consumer to use less energy. Similarly, if the utility company has aspirations to institute a demand response program, a simplified model may consider 
 \begin{equation}
    g(y)=-\left(y-\yref\right)^2
   \label{eq:demandresp}
 \end{equation}
where $\yref$ is the reference signal prescribed by the demand response program.  

The consumer's true utility is assumed to be 
\begin{align}
  J_F(\gamma(y), y)=-py+\gamma(y)+f(y)
  \label{eq:trueagent}
\end{align}
where $p$ is the price of energy set and known to all and $\gamma:Y\rar\mb{R}$ is the incentive mechanism. Thus, the consumer solves the optimization problem
\begin{equation}
  \max_{y}\{J_F(\gamma(y), y)|\ y\in Y\}.
  \label{eq:consumeropt1}
\end{equation}
We assume that the consumer is a \emph{household} who is not strategic in the sense that they take the incentive $\gamma$ and the price $p$ and optimize their utility function without strategically choosing $y$. In particular, we assume that the consumer is myopic in that he does not consider past or future incentives in his optimization problem. 

Incentives are designed by solving Problem~\ref{prob:incentive} where we assume $\Gamma$ to be the set of quadratic polynomials from $Y$ to $\mb{R}$. 

The leader does not know
the follower's satisfaction function $f(\cdot)$, and hence, must estimate it as he solves the incentive design problem. We will use the notation $\hat{f}$ for the estimate of the satisfaction function and $\hat{J}_L$ and $\hat{J}_F$ for the player's cost functions using the estimate of $f$. 

We propose an algorithm for iteratively estimating the agent's satisfaction function and choosing the incentive $\gamma(\cdot)$.  We do so by using a polynomial estimate of the agent's satisfaction function at each iteration and applying first-order optimality conditions. The use of more general sets of basis functions is left for future research.

Suppose that $\gamma^{(0)}$ and $\gamma^{(1)}$ are given \emph{a priori}. At each iteration the leader issues an incentive and observes the follower's reaction. The leader then uses the observations up to the current time along with his knowledge of the incentives he issued to estimate the follower's utility function. 

Formally, at the $k$-th iterate the leader will observe the follower's reaction $y^{(k)}$ to a delivered incentive $\gamma^{(k)}$ where we suppress the dependence of the incentive on $y$. The follower's reaction $y^{(k)}$ is optimal with respect to 
$J_F(\gamma^{(k)}(y), y)$ subject to $y^{(k)}\in Y$.

We use the observations $y^{(0)}, \ldots, y^{(k)}$ to estimate the parameters in the follower's satisfaction function given by 
\begin{equation}
  \hat{f}^{(k)}(y)=\sum_{i=0}^j\alpha_iy^{i+1}
  \label{eq:estimatesatis}
\end{equation}
where $j$ is the order of the polynomial estimate to be determined in the algorithm and we restrict $\alpha=(\alpha_1, \ldots, \alpha_j)\in \mc{A}$ a convex set, e.g. $\mc{A}=\mb{R}^{j+1}_+$.


As in \cite{keshavarz:2011aa}, we assume that an appropriate constraint qualification holds and we use Kharush-Khun-Tucker (KKT) conditions to define a notion of \emph{approximate optimality}. Thus, we can allow for some error in the estimation problem either from measurement noise or suboptimal consumer choice. In particular, for each $i=0, \ldots, k$ let $\gamma^{(i)}_{y^{(i)}}=\gamma^{(i)}(y^{(i)})$ and define
\begin{align}
  r_{\text{ineq}}^{(i)}&=(g_\ell(\gamma^{(i)}_{y^{(i)}},y^{(i)}))_+, \ \ell=1, 2\\
  r_{\text{stat}}^{(i)}(\alpha, \lambda^{(i)})&=\nabla J_F(\gamma^{(i)}_{y^{(i)}},y^{(i)})+\sum_{\ell=1}^2\lambda_\ell^{(i)}\nabla g_i(\gamma^{(i)}_{y^{(i)}}, y^{(i)})\\
  r_{\text{comp}}^{(i)}(\lambda^{(i)})&=\lambda_\ell^{(i)} g_\ell(\gamma^{(i)}_{y^{(i)}}, y^{(i)}), \ \ell=1,2
  \label{eq:resids}
\end{align}
where $(\cdot)_+=\max\{\cdot, 0\}$,
\begin{equation}
  g_1(\gamma, y)=-y\leq 0 \ \text{ and } \ g_2(\gamma, y)=y-\ymax\leq 0
  \label{eq:constraints}
\end{equation}
with Lagrange multipliers $\lambda^{(i)}=(\lambda_1^{(i)}, \lambda_2^{(i)})$.

Then, for $\{(\gamma^{(i)}, y^{(i)})\}_{i=0}^k$, we can solve 
\begin{align}\label{eq:convexopt}
  \min_{\alpha, \lambda}\Big\{ \sum_{i=0}^k \phi(r_{\text{stat}}^{(i)}, r_{\text{comp}}^{(i)})|\  \alpha \in \mc{A}, \lambda^{(i)}\geq 0, \ i=0, \ldots, k\Big\}  \tag{P-2}
  \end{align}
  where the inequality for $\lambda^{(i)}$ is element-wise and $\phi:\mb{R}\times \mb{R}^2\rar \mb{R}_+$ is a nonnegative convex penalty function (e.g. any norm on $\mb{R}\times \mb{R}^2$) satisfying
\begin{equation}
  \phi(x_1, x_2)=0\Longleftrightarrow \{x_1=0, \ x_2=0\}.
  \label{eq:phisat}
\end{equation}
The optimization problem \eqref{eq:convexopt} is convex since $r_{\text{stat}}^{(i)}$ and $r_{\text{comp}}^{(i)}$ are linear in $\alpha$ and $\lambda^{(i)}$ and the constraints are convex. If we solve \eqref{eq:convexopt} and $\phi$ is zero at the optimal solution, $r_{\text{stat}}^{(i)}$ and $r_{\text{comp}}^{(i)}$ are zero at the optimal solution for each $i$. If, in addition, $r_{\text{ineq}}^{(i)}$ is zero at the optimal solution for each $i$, then the estimate $\hat{f}^{(k)}$ at iteration $k$ is exactly consistent with the data.

If $\y{i}\in \mathring{Y}$, the problem simplifies to a checking a linear algebra condition.
Indeed, consider
\begin{equation}
  \hat{J}_F(\gamma(y), y)=-py+\gamma(y)+\hat{f}(y).
  \label{eq:estimateJa}
\end{equation}
In the case that $\hat{f}$ is concave and under our assumption that the follower is rational and hence plays optimally, the observation $y^{(i)}$ is a global optimum at iteration $i$. 
Otherwise, the observation $y^{(i)}$ is a local optimum; the follower plays myopically. In both cases, we use the necessary condition \citep{Bertsekas:1999fk}
\begin{equation}
  \nabla \hat{J}_F^{(i)}(\gamma(y^{(i)}), y^{(i)})=0  
  \label{eq:nece2}
\end{equation}
for each of the past iterates $i\in\{0, \ldots, k\}$ to determine estimates of the coefficients in $\hat{f}^{(k)}$. 

At the $k$-th iteration, we have data $\{(\gamma^{(i)}, y^{(i)})\}_{i=0}^k$. Since we require $\gamma^{(i)}\in \Gamma$, we can express each $\gamma^{(i)}$ as
\begin{equation}
  \gamma^{(i)}(y)=\xi_1^{(i)}y+\xi_2^{(i)}y^2.
  \label{eq:gammai}
\end{equation}
Then, using Equation \eqref{eq:nece2}, we define
\begin{equation}
  b^{(i)}=p-(\xi_1^{(i)}+2\xi_2^{(i)}y^{(i),d})
  \label{eq:bi}
\end{equation}
and
\begin{equation}
  \tilde{y}_{j}^{(i)}=\bmat{1& 2\y{i} & \cdots& (j+1)(\y{i})^j}
  \label{eq:ytilde}
\end{equation}
for $i\in \{0, \ldots, k\}$. 

We want to find the lowest order polynomial estimate of $f$ given the data.
We do so by checking if $\mathrm{b}^{(k)}\in \text{range}(Y^{(k)})$ where
\begin{equation}
   Y^{(k)}=\bmat{-& \tilde{y}_{j}^{(0)} &-\\ &\vdots
  &\\-& \tilde{y}_{j}^{(k)} &-},\ \mathrm{b}^{(k)}=\bmat{b^{(0)}\\ \vdots \\ b^{(k)}}
  \label{eq:range}
\end{equation}
starting with $j=2$ and increasing it until \eqref{eq:range} is satisfied or we reach $j=k$. Suppose that it is satisfied at $j=N$, $2\leq N\leq k$. Then, we estimate $\hat{f}^{(k)}$ to be an $(N+1)$-th order polynomial. We determine $\alpha_i$ for $i\in\{0, \ldots, N\}$ by solving
\begin{equation}
  \mathrm{b}^{(k)}-Y^{(k)}\alpha =0, \ \text{where }\ \alpha=[\alpha_0\ldots \alpha_N]^T.
  \label{eq:alphen}
\end{equation}
If $\mathrm{b}^{(k)}\notin \text{range}(Y^{(k)})$ for any $j\in [2,k]$, we terminate.

Our algorithm prescribes that the leader check if $y^{(i)}\in \mathring{Y}$ for each $i$. If this is the case, then he shall find the minimum order polynomial given the data as described above. On the other hand, he shall solve the convex problem \eqref{eq:convexopt}. 

Using $\{\alpha_i\}_{i=0}^j$, $\hat{J}_L^{(k)}$, and $\hat{J}_F^{(k)}$, the leader solves the incentive design problem. That is, the leader first solves
\begin{align}
(v^{(k+1),d}, y^{(k+1),d})&=\arg\min_{v\in V,y\in Y}\hat{J}_L^{(k)}(v,y)\\
&=\arg\min_{v\in V,y\in Y}\left\{g(y)-v+\beta \hat{f}^{(k)}(y)\right\}
  \label{eq:principopt}
\end{align}
Then, the leader finds $\gamma^{(k+1)}\in \Gamma$ such that
\begin{eqnarray}
  \argmax_{y\in Y} \hat{J}^{(k)}_F(\gamma^{(k+1)}(y), y)&= y^{(k+1),d}\label{eq:incentive3-2}\\
  \gamma^{(k+1)}(y^{(k+1),d})&=v^{(k+1),d}
  \label{eq:incentive3}
\end{eqnarray}
If $y^{(k+1),d}\in \mathring{Y}$, then since we restrict $\gamma^{(k+1)}$ to be of the form \eqref{eq:gammai} the above problem reduces to solving $A^{(k+1)}\xi^{(k+1)}=\tilde{\mathrm{b}}^{(k+1)}$ where
\begin{align}
  A^{(k+1)}=\begin{bmatrix} 1 & 2y^{(k+1), d}\\ y^{(k+1),d} & (y^{(k+1),d})^2\end{bmatrix}, 
    \xi^{(k+1)}=\begin{bmatrix} \xi^{(k+1)}_1\\ \xi^{(k+1)}_2\end{bmatrix}, 
  \label{eq:Axi}
\end{align}
and 
\begin{equation}
  \tilde{\mathrm{b}}^{(k+1)}=\begin{bmatrix}p-\alpha_0-2\alpha_1y^{(k+1),d}\\ v^{(k+1),d}\end{bmatrix}.
  \label{eq:btilde}
\end{equation}
If $\tilde{\mathrm{b}}^{(k+1)}\in \text{range}(A^{(k+1)})$ then a solution $\xi^{(k+1)}$ exists and if $A^{(k+1)}$ is full rank
, the solution is unique. Otherwise, if $y^{(k+1),d}\notin \mathring{Y}$,
we terminate the algorithm.
\begin{remark}
  The algorithm is motivated by the case when the consumer's satisfaction function is a polynomial of order $k$ and the utility company does not know $k$, by following the algorithm past even $k+1$ iterations, the utility company will be playing optimally. Alternatively, if incentives $\gamma^{(i)}$ were chosen randomly, the utility company would not know when to stop choosing random $\gamma^{(i)}$'s; thus, after $k+1$ iterations would begin playing suboptimally. 
\end{remark}

  \begin{proposition}
    Let $f$ be polynomial of order $k+1$, $y^{(0)}\in \mathring{Y}$ and $\gamma^{(0)}, \gamma^{(1)}$ be given \emph{a priori}. Suppose that at each iteration of the algorithm $\mathrm{b}^{(\ell)}\in \text{range}(Y^{(\ell)})$, $\text{rank}(Y^{(\ell)})=\ell+1$, $y^{(i)}\in \mathring{Y}$  and $\tilde{\mathrm{b}}^{(\ell+1)}\in \text{range}(A^{(\ell+1)})$.  Then, after $k$ iterations 
 the satisfaction function is known exactly and the incentive $\gamma^{(k+1)}$
 induces the consumer to use the desired control.
\end{proposition}

\begin{proposition}
  Suppose that $f$ is polynomial up to order $k+1$ and that the leader has $k+1$ historical measurements
  \begin{equation}
    \gamma^{(-k)}, \ldots, \gamma^{(1)}, y^{(-k)}, \ldots, y^{(1)}
    \label{eq:historical}
  \end{equation}
  such that $Y^{(k)}$ is full rank where $y^{(i)}\in \mathring{Y}$ for $i=0, \ldots, k$, then the leader can estimate the follower's satisfaction function exactly and if there exists an incentive $\gamma^{(k+1)}$, then it induces the desired equilibrium.
 \end{proposition}

We conclude this section by providing an example of the iterative process when $f$ is a concave function.
  \begin{example}
    First, we suppose that $\gamma^{(0)},\gamma^{(1)}\in \Gamma$ are chosen \emph{a priori} and are parameterized as follows:
    \begin{align}
      \gamma^{(i)}(y)=\xi_{1}^{(i)}y+\xi_2^{(i)}y^2
      \label{eq:gammazerone}
    \end{align}
    for $i\in\{0,1\}$.
    Then, the procedure goes as follows. The leader issues $\gamma^{(0)}$ and observes $y^{(0)}$. Subsequently, he issues $\gamma^{(1)}$ and observes $y^{(1)}$. Suppose $y^{(0)}, y^{(1)}\in \mathring{Y}$. The leader determines $\alpha_1$, $\alpha_0$ in the estimation of 
$\hat{f}(y)=\alpha_1y^2+\alpha_0y$ by computing the derivative of $\hat{J}_F^{(0)}(\gamma^{(0)}(y),y)$ and $\hat{J}_F^{(1)}(\gamma^{(1)}(y),y)$
with respect to $y$, evaluating at $\y{0}$ and $\y{1}$ and equating to zero, i.e. he solves
\begin{eqnarray}
  -p-2\y{0}+2(\alpha_1+\xi_2^{(0)})\y{0}+\alpha_0+\xi_1^{(0)}&=& 0\\
  -p-2\y{1}+2(\alpha_1+\xi_2^{(1)})\y{1}+\alpha_0+\xi_1^{(1)}&=& 0
  \label{eq:eqns}
\end{eqnarray}
for $\alpha_0$ and $\alpha_1$. If either $y^{(0)}$ or $y^{(1)}$ are on the boundary of $Y$, then the leader solves \eqref{eq:convexopt} for $\alpha=(\alpha_1, \alpha_0)$. 

Using $\alpha_0, \alpha_1$, the leader solves the following incentive design problem for $\gamma^{(2)}$. First, find $(v^{(2),d}, y^{(2),d})\in V\times Y$ such that
\begin{equation}
  \hat{J}_L^{(2)}(v,y)=-y-v+\alpha_1y^2+\alpha_0y
  \label{eq:Jptwo}
\end{equation}
 is maximized. Since we restrict to quadratic incentives, we parameterize $\gamma^{(2)}$ as in Equation~\eqref{eq:gammazerone} with $i=2$. 
Now, given the utility $\hat{J}_F^{(2)}(\gamma^{(2)}(y), y)$,
we find $\xi_1^{(2)}, \xi_2^{(2)}$ such that 
\begin{eqnarray}
  \arg\max_{y\in Y}\hat{J}_F^{(2)}(y;\xi_1^{(2)}, \xi_2^{(2)} )&=& y^{(2),d}\label{eq:argmaxJ}\\
  \xi_1^{(2)}y^{(2),d}+\xi_2^{(2)}(y^{(2),d})^2&=& v^{(2),d}
  \label{eq:incentivetwo}
\end{eqnarray}
Assuming that $y^{(2),d}\in \mathring{Y}$, it will be
an induced local maxima under the incentive $\gamma^{(2)}$. Hence, Equation~\eqref{eq:argmaxJ} can be reformulated using the necessary condition 
\begin{equation}
  \nabla_y\hat{J}_F^{(2)}(y^{(2),d}; \xi_1^{(2)}, \xi_2^{(2)} )=0.
  \label{eq:nece3}
\end{equation}
Now, Equations~\eqref{eq:incentivetwo} and \eqref{eq:nece3} give us two equations in the two unknowns $\xi_1^{(2)}, \xi_2^{(2)}$ that can be solved; indeed,
\begin{align}
  -p+\xi_1^{(2)}+\alpha_0+2(\xi_2^{(2)}+\alpha_1)y^{(2),d}&=0\\
  \xi_1^{(2)}y^{(2),d}+\xi_2^{(2)}(y^{(2),d})^2&= v^{(2),d}
  \label{eq:twoeqns}
\end{align}
Solving these equations gives us the parameters for $\gamma^{(2)}$. Now, the leader can issue $\gamma^{(2)}$ to the follower and observe his reaction $y^{(2)}$. The leader can then continue in the iterative process as described above.
  \end{example}

\section{Device Level Incentive Design Using Disaggregation Algorithm}
In a manner similar to the previous section, we consider that the consumer's satisfaction function is unknown. However, we now consider that the utility company desires to design device level incentives. We remark that the utility company may not want to incentivize every device; the process we present can be used to target devices with the highest consumption or potential to offset inaccuracies in load forecasting.

\subsection{Exact Disaggregation Algorithm}
\label{sec:disaggregation_a}
We first describe the process of designing device level incentives assuming the utility company has a disaggregation algorithm in place which produces no error. That is, they observe the aggregate signal and then applies their disaggregation algorithm to get exact estimates of the device level usage $y_\ell$ for $\ell\in \{1, \ldots, D\}$ where $D$ is the number of devices. 

The utility company has the true utility function
\begin{align}
  J_L(v,y)=\sum_{\ell=1}^Dg_\ell(y_\ell)-v_\ell+\beta_\ell f_\ell(y_\ell)
  \label{eq:trueprincipdl}
\end{align}
and the consumer has the true utility function
\begin{align}
  J_F(\gamma(y), y)=\sum_{\ell=1}^D-py_\ell+\gamma_\ell(y_\ell)+f_\ell(y_\ell).
  \label{eq:trueagentdl}
\end{align}
The utility company could choose only to incentivize specific devices such as high consumption devices. This fits easily into our framework; however, for simplicity we just present the model in which incentives are designed for each device. 

The implicit assumption that the player utilities are separable in the devices allows us to generalize the algorithm presented in the previous section. Let us be more precise. We again assume that $\gamma_\ell^{(0)}, \gamma_\ell^{(1)}$ for $\ell\in \{1, \ldots, D\}$ are given \emph{a priori}. 

At the $k$-th iteration the utility company issues an incentive $\gamma_\ell^{(k)}$ for each device $\ell\in\{1, \ldots, D\}$ and observes the aggregate signal $y^{(k)}$. Then they apply a disaggregation algorithm to determine the device level usage $y_\ell^{(k)}$ for $\ell\in\{1, \ldots, D\}$.

The utility company forms an estimate of the consumer's device level satisfaction function
\begin{equation}
  \hat{f}^{(k)}_\ell(y_\ell)=\sum_{i=0}^j\alpha_{i,\ell}y_\ell^{i+1}
  \label{eq:fhatkdl}
\end{equation}
and then solves the problem of finding the $\alpha_{i,\ell}$'s by solving for $\alpha_\ell=(\alpha_{0,\ell}, \ldots, \alpha_{j, \ell})$ as in the previous section for each device $\ell\in \{1, \ldots, D\}$.

\begin{proposition}
  For $\ell\in\{1, \ldots, D\}$, let $f_\ell$ be polynomial up to order $k_\ell+1$, $\gamma^{(0)}_\ell, \gamma^{(1)}_\ell$ be given \emph{a priori}, and $y^{(0)}_\ell, y^{(1)}_\ell\in \mathring{Y}$. Suppose that at each iteration of the algorithm $\mathrm{b}^{(m)}_\ell\in \text{range}(Y_\ell^{(m)})$, $\text{rank}(Y_\ell^{(m)})=m_\ell+1$, $y^{(m)}_\ell\in \mathring{Y}$, and $\tilde{\mathrm{b}}^{(m+1)}_\ell\in A^{(m+1)}_\ell$ for each $\ell\in\{1, \ldots, D\}$. 
Then, after 
\begin{equation}
  k^\ast=\max_{\ell\in \{1, \ldots, D\}}k_\ell
  \label{eq:kstar}
\end{equation}
iterations, the satisfaction function is known exactly and the incentives $\gamma^{(k^\ast+1)}_\ell$ induce the desired equilibrium.
\end{proposition}
Note that the notation $(\cdot)_\ell$ indicates the object defined in Section~\ref{sec:aggregate} for the $\ell$-th device. 

\subsection{Disaggregation Algorithm with Some Error}
\label{sec:disaggregation_b}
Now, we consider that the leader has some error in his estimate of the device level usage due to inaccuracies in the disaggregation algorithm, i.e. the leader determines $\hat{y}_\ell$ such that $\|y_\ell-\hat{y}_\ell\|\leq \vep$ where $\vep>0$ is the resulting error from the estimation in the disaggregation algorithm. Bounds on $\vep$ can be determined by examining the fundamental limits of non--intrusive load monitoring algorithms \citep{dong:2013aa}.


We again assume that $\gamma_\ell^{(0)}, \gamma_\ell^{(1)}$ for $\ell\in \{1, \ldots, D\}$ are given \emph{a priori}. Following the same procedure as before, at the $k$-th iterate the leader will issue $\gamma_\ell^{(k)}$ for $\ell\in\{1, \ldots, D\}$ and observe $y^{(k)}$. Then apply a disaggregation algorithm to determine $\hat{y}_\ell^{(0)}$ where 
\begin{equation}
  \|y_\ell-\hat{y}_\ell\|\leq \vep
  \label{eq:estimate}
\end{equation}
for $\ell\in\{1, \ldots, D\}$. The incentive design problem follows the same steps as provided in the previous section
with the exception that the $y_\ell$'s are replaced with the estimated $\hat{y}_\ell$'s and we tolerate an error in solving 
for the minimal polynomial estimate of $f_\ell$.
\begin{figure}[ht]
  \begin{center}
    \includegraphics[width=\columnwidth]{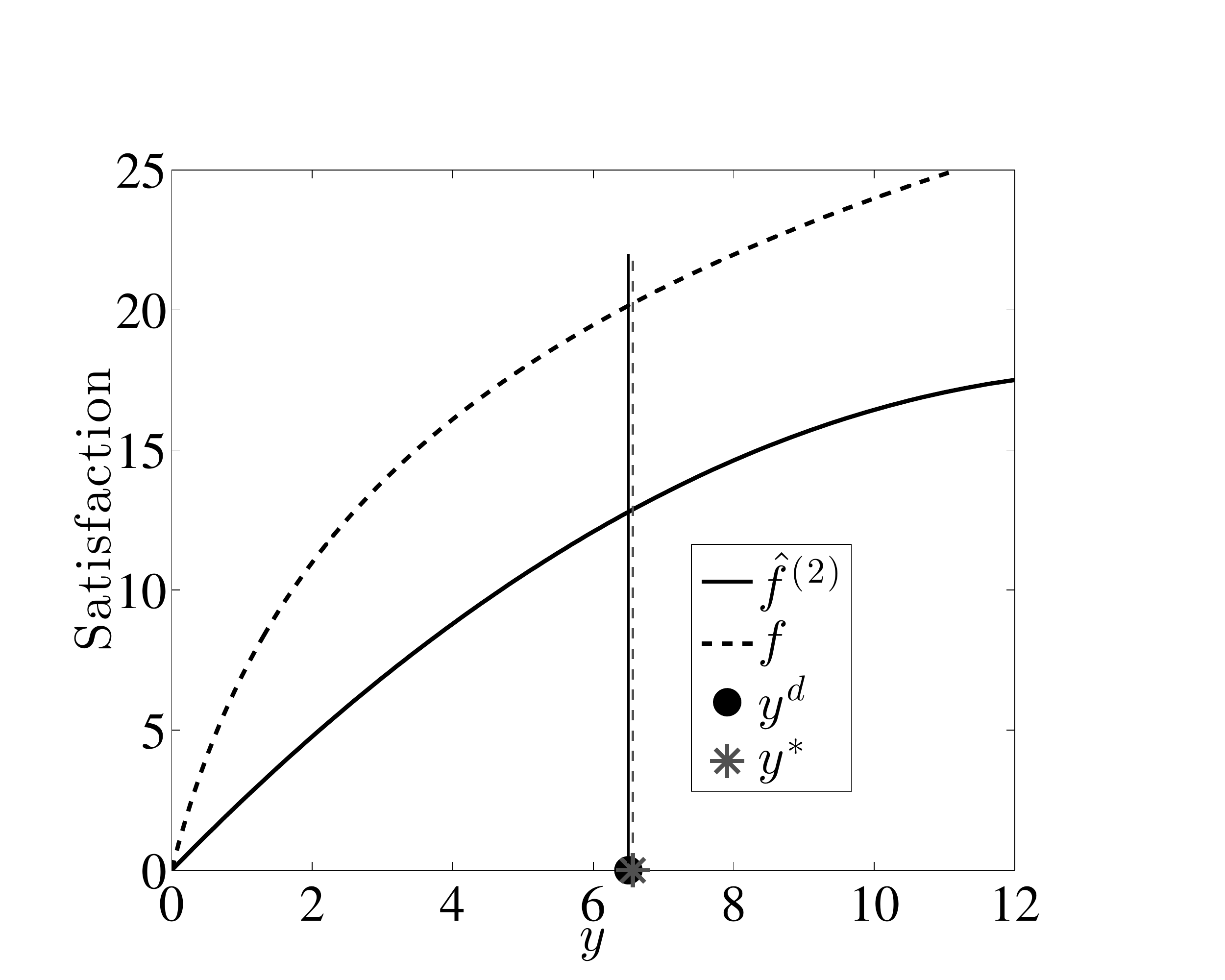}
  \end{center}
  \caption{Estimated satisfaction function $\hat{f}^{(2)}$ and true satisfaction function $f$. The true response $y^\ast=6.56$ and the desired response $y^d=6.5$. Notice that the slope of the estimated satisfaction function and the slope of the true satisfaction function are roughly equal at $y^d$ and $y^\ast$.}
  \label{fig:satisfaction}
\end{figure}

\section{Numerical Examples}
\label{sec:numerical}
We simulate two examples of designing incentives while estimating the consumer's satisfaction function. 
In both examples we assume a unit price per unit of energy, i.e. $p=1$.

\subsection{Aggregate Signal and Log Satisfaction Function}
We simulate a system in which the consumer has the true utility given by
\begin{equation}
  J_F(\gamma(y), y)=-py+\gamma(y)+f(y)
  \label{eq:trueJaex1}
\end{equation}
where the satisfaction function is
\begin{equation}
  f(y)=10\log(y+1).
  \label{eq:satisfactionfn}
\end{equation}
We assume the utility company is in a regulated market; hence wants the consumer to consume less.  Thus, the utility company has utility function 
\begin{equation}
  J_L(v,y)=-y-v+\beta f(y)
  \label{eq:jpex1}
\end{equation}
where the benevolence factor is $\beta=0.75$. We let $\ymax=\vmax=100$.
We choose two concave incentive function $\gamma^{(0)}(y)$ and $\gamma^{(1)}(y)$ defined as follows: $\gamma^{(0)}(y)=-y^2+10y$ and $\gamma^{(1)}(y)=-y^2+15y$.

We use the algorithm presented in Section~\ref{sec:aggregate} to design incentives while estimating $\alpha_0$ and $\alpha_1$. We simulate the utility company issuing $\gamma^{(0)}$ and then $\gamma^{(1)}$ where the consumer chooses his optimal response to each of the incentives. The responses are $y^{(0)}=5.29$ and $y^{(1)}=7.58$.
After two iterations, we get a reasonable approximation of the true $f$ and a quadratic incentive $\gamma^{(2)}$;
\begin{equation}
  \hat{f}^{(2)}(y)=2.57y-0.093y^2, \ \gamma^{(2)}(y)=0.33y-0.05y^2.
  \label{eq:fhat2log}
\end{equation}

The optimal power usage under the incentive $\gamma^{(2)}$ is $y^\ast=6.56$ and the desired power usage is $y^d=6.5$. It is clear that the utility company could do better if he new the true satisfaction function. Figure~\ref{fig:satisfaction} shows $\hat{f}^{(2)}(y)$ and $f(y)$. It is important to notice that $y^\ast$ is nearly equal to $y^d$ and at these two points the slope of $\hat{f}^{(2)}$ is approximately equal to that of the true $f$. This indicates that $\hat{f}^{(2)}$ is a good estimate of $f$.

\begin{figure}[h]
  \begin{center}
    \vspace{-0.1in}
    \includegraphics[width=\columnwidth]{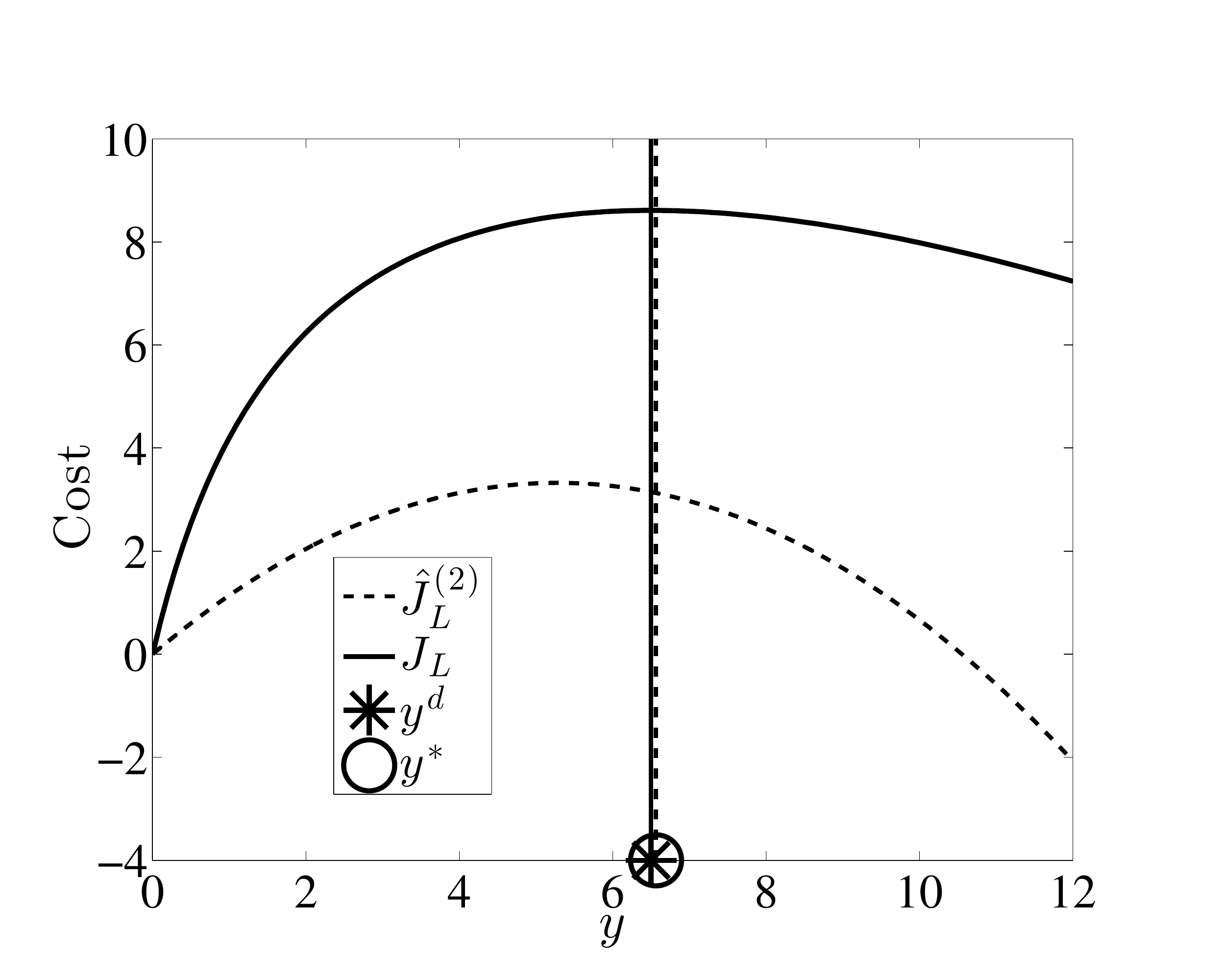}
  \end{center}
  \caption{Estimated cost function $\hat{J}_L^{(2)}(\gamma^{(2)}(y), y)$ and true cost function $J_L(v^d,y)$ along $v=v^d$ for simulation with $\log$--satisfaction function.
  }
  \label{fig:costs}
\end{figure}

Figure~\ref{fig:costs} shows the true utility function of the utility company $J_L(v^d, y)$ with $v=v^d$ fixed and the estimated utility $\hat{J}_L^{(2)}(\gamma^{(2)}(y), y)$. $y^d$ is the point at which $J_L(v^d, y)$ is maximized and it is approximately the point where $\hat{J}_L^{(2)}(\gamma^{(2)}(y), y)$ is maximized. It is important to note the \emph{shape} of $J_L$ and $\hat{J}_L^{(2)}$. The offset is not important because we are not estimating a constant term in $\hat{f}$ since it does not affect the optimal response, i.e. if you shift $\hat{J}_L^{(2)}$ by a constant term, $y^\ast$ is still the optimal response. 

\subsection{Disaggregated Signal}
We simulate a system in which the consumer's the true utility 
\begin{equation}
  J_F(\gamma(y), y)=\sum_{\ell=1}^{10}-py_\ell+\gamma_\ell(y_\ell)+f_\ell(y_\ell)
  \label{eq:trueJadis}
\end{equation}
where the satisfaction functions $f_\ell(y_\ell)$ are exactly quadratic for each device $\ell\in\{1,\ldots, 10\}$;
\begin{equation}
  f_\ell(y_\ell)=\alpha_{1,\ell}y_\ell^2+\alpha_{0,\ell} y_\ell
  \label{eq:quadsatis}
\end{equation}
The utility company's utility is given by
\begin{equation}
  J_L(v,y)=\sum_{\ell=1}^{10}-y_\ell-v_\ell+\beta_\ell f_\ell(y_\ell)
  \label{eq:jptruedis}
\end{equation}
where the benevolence factor (i.e. a representation of how much the utility company cares about the satisfaction of the consumer) is $\beta_\ell=1$ for each $\ell$.
The utility company must disaggregate the aggregated energy signal $y$ giving rise to estimates $\hat{y}_\ell$. If $\hat{y}_\ell=y_\ell$, i.e. there is no error in the disaggregation algorithm, then after two iterations the utility company would know the satisfaction function of each device exactly. Let's explore the case when the disaggregation algorithm has $\vep$--error. In our examples we randomly generate noise within a given $\vep$ bound and add that to the true $y_i$'s to simulate the error in the disaggregation step resulting from the disaggregation algorithm. 
\begin{figure}[ht]
  \begin{center}
    \includegraphics[width=\columnwidth]{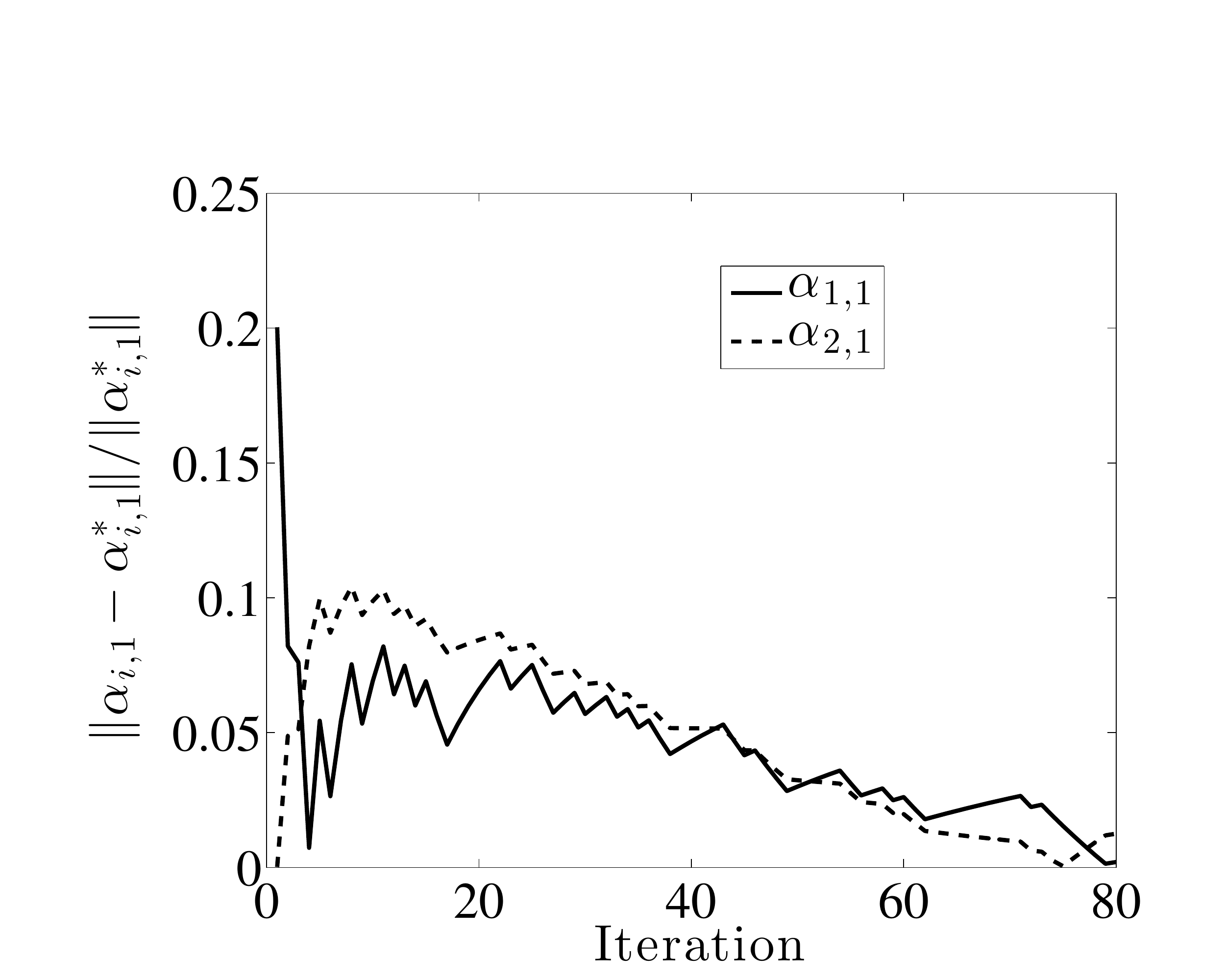}
  \end{center}
  \caption{Relative error in estimate of $\alpha_{i,1}$'s for device 1 with disaggregation error bound $\vep=0.15$. $\alpha_i^\ast$ is the true value.  The relative error eventually decreases below the noise bound $\vep=0.15$.}
  \label{fig:error2}
\end{figure}

Figure~\ref{fig:error2} shows the relative error on the estimates of $\alpha_{i,1}$ for $i\in \{1,2\}$ as a function of the iteration. 
The relative error for other devices are similar. We used the error bound $\vep=0.15$ for the disaggregation error. The relative error decreases as the number of iterations increase.
It eventually ends up below the noise bound $\vep$ and remains there. 

As we iterate the noise introduced via disaggregation has minimal impact on the estimate of $\alpha_{i, \ell}$ for $i=\{1, 2\}$ and $\ell\in\{1, \ldots, D\}$. We note that the designed incentive for this problem converges to zero as we increase the iterations and the impact of the noise is minimized. 
It becomes zero since the benevolence factor is $\beta_\ell=1$ and the price $p=1$; hence, the agent and the leader have the same utility functions after the leader learns the agent's satisfaction function. As we increase the noise threshold $\vep$, the estimation of $\alpha_{i, \ell}$ degrades. 

%
In the last simulation, we decrease the benevolence factor to $\beta=0.75$ and comment on the resulting incentives.
\begin{figure}[ht]
  \begin{center}
    \includegraphics[width=\columnwidth]{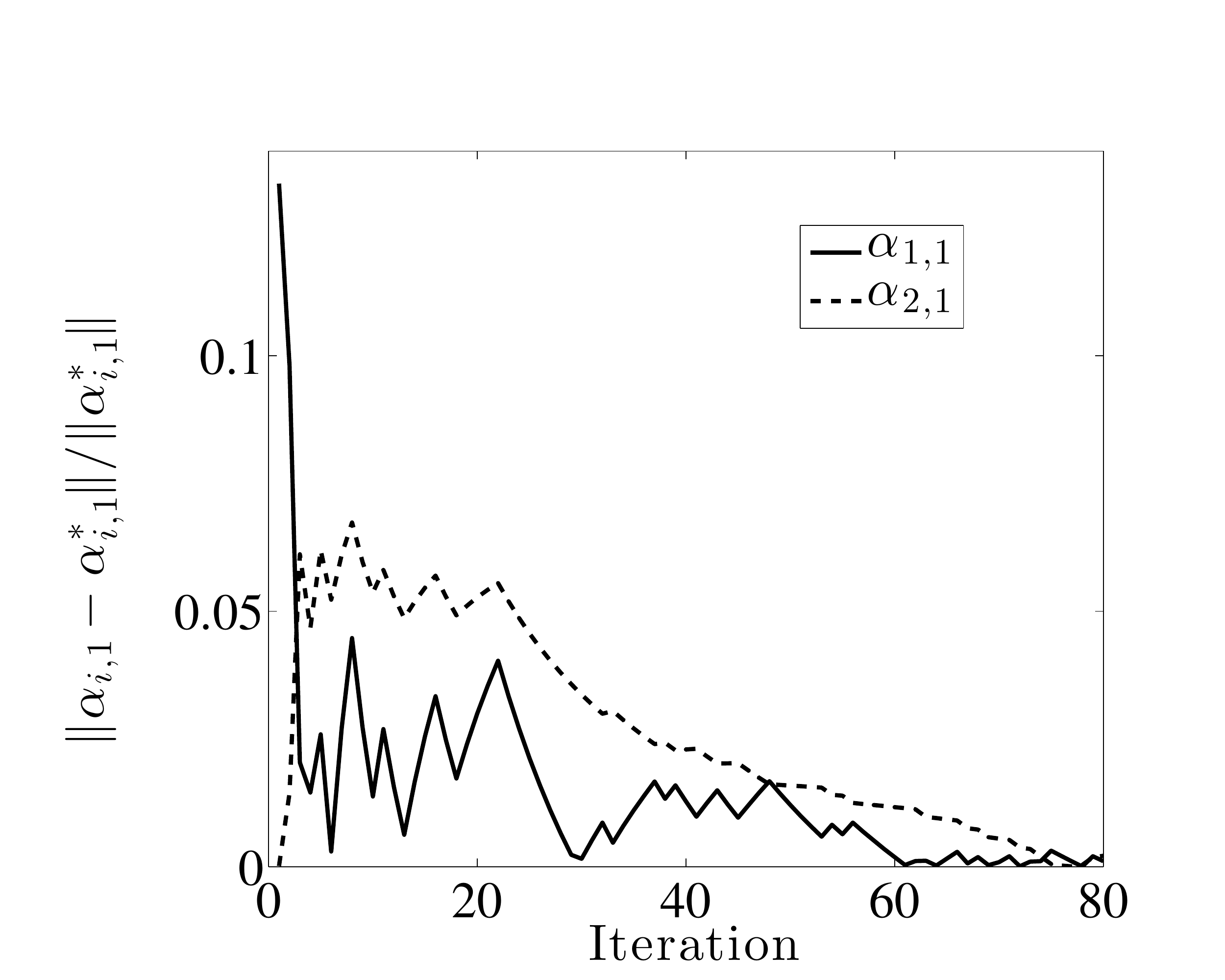}
  \end{center}
  \caption{Relative error in estimate of $\alpha_{i,1}$'s for device 1 with disaggregation error bound $\vep=0.1$. $\alpha_i^\ast$ is the true value.  The relative error eventually decreases below the noise bound $\vep=0.1$;}
  \label{fig:error1}
\end{figure}

In Figure~\ref{fig:error1} we show again the relative error in the estimates of $\alpha_{i, 1}$ for $i\in\{1, 2\}$. The relative error of each of the other devices is similar. In this case the incentive for device $\ell=1$ converge to
\begin{equation}
  \gamma^{\ast}_1(y)=-0.39y^2+0.33y.
  \label{eq:incentive1}
\end{equation}
The other devices have similar incentives. 
The reason that there is a non--zero incentive is due to the fact that the utility company is not completely benevolent; he does not care as much about the satisfaction of the consumer as he does the other terms in his cost function. However, as we iterate and the utility company learns the consumer's cost function, the utility company is able to use the incentives to force the desired action $y^\ast=y^d$ where $y^\ast$ is the consumer's true response.
\section{Discussion and Future Work}
\label{sec:discussion}
We modeled the utility company--consumer interaction as a reversed Stackelberg game.
We defined a process by which the utility company can jointly estimate the consumer's utility function and design incentives for behavior modification. 
Whether the utility company is interested in inducing energy efficient behavior or creating an incentive compatible demand response program, the procedure we present applies. 
We are studying fundamental limits of \emph{non--intrusive load monitoring} in order to determine precise bounds on the payoff to the utility company when a disaggregation algorithm is in place and incentives are being designed. 
We seek to understand how these fundamental limits impact the quality of the incentive design problem as well as how they can be integrated into a stochastic contorl framework for incentive design when faced with non-strategic agents with unknown preferences.


The electrical grid is a social cyber-physical system (S-CPS) with human actors influencing the trajectory of the system. 
Inherent to the study of S-CPS's are privacy and security considerations.
We remark that consumers may consider their satisfaction function to be private information. We are currently exploring the design of privacy--aware mechanisms for $\vep$--incentive compatible problems \citep{nissim:2012aa}. 

\bibliography{2014ifacv2}  
\end{document}